\documentclass[12pt]{article}
\usepackage{amsmath}
\usepackage{amssymb}
\usepackage[ctr1,eng]{myarticl}

\newcommand{\Aut }{\mathrm{Aut}}

\newcommand{\cB }{\mathcal{B}}

\newcommand{\End }{\mathrm{End}}

\newcommand{\lact }{.}

\newcommand{\lcoa }{\delta }

\newcommand{\Ndbasis }{\boldsymbol{\mathrm{e}}}

\newcommand{\ndN }{\mathbb{N}}

\newcommand{\ndZ }{\mathbb{Z}}

\newcommand{\ot }{\otimes }

\newcommand{\PBW }{Poincar\'e--Birkhoff--Witt }
\newcommand{\proots }{\boldsymbol{\Delta }^+}

\newcommand{\qnum }[2]{(#1)_{#2}}

\newcommand{\roots }{\boldsymbol{\Delta }}

\newcommand{\YD }{Yetter--Drinfel'd }

\title{Weyl equivalence for rank 2 Nichols algebras\\
of diagonal type
\thanks{Supported by the European Community under a Marie Curie
Intra-European Fellowship}}
\author{I.~Heckenberger}

\begin{document}

\maketitle

\begin{abstract}
\YD modules of diagonal type admit an equivalence relation which
preserves dimension and Gel'fand--Kirillov dimension of the
corresponding Nichols algebras. This relation is determined
explicitly for all rank 2 \YD modules where the Gel'fand--Kirillov
dimension is known to be finite.


\end{abstract}

\section{Introduction}

The method of Andruskiewitsch and Schneider \cite{a-AndrSchn98}
for the classification of pointed Hopf algebras is based in an
essential way on the knowledge of Nichols algebras of group type.
There exist several methods to determine such Nichols algebras,
see for example
\cite{a-Rosso98},\cite{a-AndrGr99},\cite{a-Khar99},\cite{a-AndrSchn00},\cite{inp-MilSchn00},\cite{a-Grana00},\cite{a-Heck04b},\cite{a-Heck04c}.
Nichols algebras of Cartan type can be described with help of root
systems. However in the slightly more general setting, where the
underlying \YD module is of diagonal type, there exists still no
satisfactory classification result.

In \cite{a-Heck04a} and \cite{a-Heck04b} the author classified
finite dimensional rank 2 Nichols algebras of diagonal type.
Compared to the theory of root systems the set of solutions is
much more variegated. This causes difficulties for a
generalization in higher rank. However, as for root systems the
Weyl group plays an essential role, \YD modules of diagonal type
admit the definition of the so called Weyl--Brandt groupoid
\cite{a-Heck04c}. The transformations corresponding to the
elements of this groupoid are usually not algebra homomorphisms
between Nichols algebras, but they preserve their dimension and
Gel'fand--Kirillov dimension. Thus one can reduce the
classification problem to the question of finding the orbits of
all possible (say finite) Weyl--Brandt groupoids. In the present
paper the equivalence classes of all known rank 2 \YD modules of
diagonal type obtained from the action of the Weyl--Brandt
groupoid are determined. Conjecturally all rank 2 \YD modules of
diagonal type, such that the corresponding Nichols algebra has
finite Gel'fand--Kirillov dimension, are twist equivalent to one
in Proposition \ref{s-rank2}. Therefore in the classification of
Nichols algebras of arbitrary rank Proposition \ref{s-rank2} can
be used to pick a convenient representant of an orbit of the
Weyl--Brandt groupoid.

If not stated otherwise the notation in \cite{a-AndrSchn00} and
\cite{a-Heck04a,a-Heck04b,a-Heck04c} is followed.

\section{The Weyl--Brandt groupoid of a Nichols algebra of
diagonal type} \label{sec-WBg}

Let $k$ be a field of characteristic zero, $G$ an abelian group,
and $V$ a \YD module over $kG$ of rank $n$ for some $n\in \ndN $.
More precisely, let $\{x_1,x_2,\ldots ,x_n\}$ be a basis of $V$,
$\{g_i\,|\,1\le i\le n\}$ a subset of $G$, and $q_{ij}\in
k\setminus \{0\}$ for all $i,j\in \{1,2,\ldots ,n\}$, such that
\begin{align*}
\lcoa (x_i)=&g_i\ot x_i,& g_i\lact x_j=&q_{ij}x_j
\end{align*}
for $i,j\in \{1,2,\ldots ,n\}$. Then $V$ is a braided vector space
\cite[Def.\,5.4]{inp-Andr02} with braiding $\sigma \in \End (V\ot
V)$, where
\begin{align*}
\sigma (x_i\ot x_j)=q_{ij}x_j\ot x_i
\end{align*}
for all $i,j\in \{1,2,\ldots ,n\}$. The Nichols algebra $\cB (V)$
is of diagonal type \cite[Def.\,5.8]{inp-Andr02} and has a $\ndZ
^n$-grading defined by
\begin{align*}
\deg x_i:=\Ndbasis _i\quad \text{for all $i\in \{1,2,\ldots ,n\}$}
\end{align*}
where $E_0:=\{\Ndbasis _1,\ldots ,\Ndbasis _n\}$ is the standard
basis of the $\ndZ $-module $\ndZ ^n$. Let $\chi :\ndZ ^n\times
\ndZ ^n\to k\setminus \{0\}$ denote the bicharacter defined by
\begin{align*}
\chi (\Ndbasis _i,\Ndbasis _j):=q_{ij},
\end{align*}
where $i,j\in \{1,2,\ldots ,n\}$. Kharchenko proved \cite[Theorem
2]{a-Khar99} that the algebra $\cB (V)$ has a (restricted) \PBW
basis consisting of iterated skew-commutators of the elements
$x_i$ of $V$. Moreover, one can assume that this basis has the
additional property that
\begin{itemize}
\item[(P)]
the height of a \PBW generator of $\ndZ ^n$-degree $d$ is finite
if and only if $2\le \mathrm{ord}\,\chi (d,d)<\infty $
($\mathrm{ord}$ means order with respect to multiplication), and
in this case it coincides with $\mathrm{ord}\,\chi (d,d)$.
\end{itemize}
Indeed, if $\chi (d,d)$ is a primitive $m^\mathrm{th}$ root of
unity for some $m\ge 2$, but the height of the element of degree
$d$ is infinite, then one can add the $m^\mathrm{th}$ power of
this element to the \PBW basis. Its height is automatically
infinite since $\chi (md,md)=\chi (d,d)^{m^2}=1$.

Let $\proots (\cB (V))$ denote the set of degrees of the
(restricted) \PBW generators counted with multiplicities. Note
that this definition is independent of the choice of a $\ndZ
^n$-graded \PBW basis satisfying property (P).

Assume that $i\in \{1,2,\ldots ,n\}$ such that for all $j\in
\{1,2,\ldots ,n\}$ with $j\not=i$ the numbers
\begin{align}\label{eq-mij}
 m_{ij}:=&
 \min\{m\in \ndN _0\,|\,\qnum{m+1}{q_{ii}}(q_{ii}^m
 q_{ij}q_{ji}-1)=0\}
\end{align}
are well-defined. In this case one can introduce a $\ndZ $-linear
mapping $s_i:\ndZ ^n\to \ndZ ^n$ as follows:
\begin{align}\label{eq-si}
s_i(\Ndbasis _j):=&
 \begin{cases}
 -\Ndbasis _i & \text{if $j=i$,}\\
 \Ndbasis _j+m_{ij}\Ndbasis _i & \text{if $j\not=i$.}
 \end{cases}
\end{align}
The map $s_i$ is a pseudo-reflection
\cite[Ch.~5,\S2]{b-BourLie4-6}, that is $\mathrm{rk}\,(s_i-\id
)=1$, and it satisfies the equation $s_i^2=\id $.

\begin{defin}
 Let $E$ be a subset of $\ndZ ^n$. We say that a \YD module
$V'$ \textit{has degree} $E$ \textit{with respect to} $V$ if there
exists a $\ndZ ^n$-graded \YD submodule $V''$ of $\cB (V)\ot \cB
(V^*)$ isomorphic to $V'$ and a basis $\{v''_e\,|\,e\in E\}$ of
$V''$ indexed by $E$ such that $\deg v''_e=e$ for all $e\in E$. In
this case a subset $\{v'_e\,|\,e\in E\}$ of $V'$ is said to be a
\textit{basis with respect to} $V$ \textit{and} $E$ if the linear
map from $V''$ to $V'$, defined by $v''_e\mapsto v'_e$ for all
$e\in E$, is an isomorphism of \YD modules.
\end{defin}

Note that $E$ is not uniquely determined by $V'$ but for given $E$
and $V$ there is up to isomorphism a unique \YD module which has
degree $E$ with respect to $V$.

Let $E$ be an ordered basis of $\ndZ ^n$ and let $V'$ be a \YD
module of degree $E$ with respect to $V$. Choose $i\in
\{1,2,\ldots ,n\}$ and an ordered basis $B'$ of $V'$ with respect
to $V$ and $E$. Define $s_{i,E}$ (if it exists) to be the map
$s_i\in \Aut (\ndZ ^n)$ corresponding to $\cB (V')$ with respect
to the basis $B'$ of $V'$. Note that if $E=E_0$ then $s_{i,E}=s_i$
as defined above. Further, by \cite[Proposition 1]{a-Heck04c} one
also has $s_{i,s_i(E_0)}=s_i$.

In \cite{a-Heck04c} the Weyl--Brandt groupoid $W(V)$ associated to
$\cB (V)$ was defined as follows.

\begin{defin}
Let $E_0:=(\Ndbasis _1,\ldots ,\Ndbasis _n)$ denote the ordered
standard basis of $\ndZ ^n$. Define
\begin{align*}
W(V)&:=\{(s,E)\,|\,\text{$s\in \Aut (\ndZ ^n)$, $E$ is an ordered
basis of $\ndZ ^n$,}\\
&\text{there exist $m_1,m_2\in \ndN _0$, $m_1\le m_2$},
\text{$i_1,\ldots ,i_{m_2}\in \{1,2,\ldots ,n\}$,}\\
& \text{and for $1\le m\le m_2$ ordered bases $E_m$ of $\ndZ ^n$,
such that}\\
&\text{$s_{i_{m+1},E_m}$ is well defined and
$s_{i_{m+1},E_m}(E_m)=E_{m+1}$ for all $m<m_2$;}\\
&\text{$s=s_{i_{m_2},E_{m_2-1}}\cdots s_{i_{m_1+1},E_{m_1}}$,
$E=E_{m_1}$\}}.
\end{align*}
The set $W(V)$ admits a partial binary operation $\circ $ given by
\begin{align*}
(s,E)\circ (t,F):=(st,F)
\end{align*}
whenever $(s,E),(t,F)\in W(V)$ and $t(F)=E$. Then $(W(V),\circ )$
is a Brandt groupoid, see \cite[Sect. 3.3]{b-ClifPres61} or
\cite[Sect.~5]{a-Heck04c}, and is called the \textit{Weyl--Brandt
groupoid of} $\cB (V)$.
\end{defin}

There is a natural action of the Weyl--Brandt groupoid $W(V)$ on
the set
\begin{align*}
E(V):=\{E\,|\,&\text{$E$ is an ordered basis of $\ndZ ^n$,}\\
&\text{$(s,E)\in W(V)$ for some $s\in \Aut (\ndZ ^n)$}\}
\end{align*}
given by the formula
\begin{align*}
(s,E)(F):=&\begin{cases}
 s(E) & \text{if $E=F$,}\\
 \text{not defined} & \text{otherwise}
\end{cases}
\end{align*}
for all $(s,E)\in W(V)$ and $F\in E(V)$.

Let $V'$ be a \YD module of diagonal type which has degree
$s_i(E_0)$ with respect to $V$, and let $\{x'_1,\ldots ,x'_n\}$ be
a basis of $V'$ with respect to $V$ and $s_i(E_0)$. Then the
structure constants $q'_{jl}$ of $V'$, where $j,l\in \{1,2,\ldots
,n\}$, depend on the numbers $q_{jl}$ as follows:
\begin{align*}
q'_{jl}=q_{ii}^{m_{ij}m_{il}}q_{il}^{m_{ij}}q_{ji}^{m_{il}}q_{jl}
\end{align*}
for all $j,l\in \{1,2,\ldots ,n\}$, where $m_{ii}=-2$. For
$j\not=i$ set
\begin{align*}
p_{ij}:=\begin{cases}
 1 & \text{if $q_{ii}^{m_{ij}}q_{ij}q_{ji}=1$,}\\
 q_{ii}^{-1}q_{ij}q_{ji} & \text{otherwise.}
\end{cases}
\end{align*}
Then from the Equations (\ref{eq-mij}) one obtains that
\begin{align}\label{eq-qtrafo}
q'_{ii}=&q_{ii},& q'_{jj}=&p_{ij}^{m_{ij}}q_{jj},&
q'_{ij}q'_{ji}=&p_{ij}^{-2}q_{ij}q_{ji},&
q'_{jl}q'_{lj}=&p_{ij}^{m_{il}}p_{il}^{m_{ij}}q_{jl}q_{lj}
\end{align}
for all $j,l\in \{1,2,\ldots ,n\}\setminus \{i\}$. These formulas
allow the following definitions.

\begin{defin}\label{d-Weq}
Let $V',V''$ be two rank $n$ braided vector spaces of diagonal
type. If with respect to certain bases their structure constants
$q'_{jl}$ and $q''_{jl}$, where $j,l\in \{1,2,\ldots ,n\}$,
satisfy the equations
\begin{align}
q'_{jj}=&q''_{jj},& q'_{jl}q'_{lj}=&q''_{jl}q''_{lj}
\end{align}
for all $j,l\in \{1,2,\ldots ,n\}$ then $W(V')=W(V'')$. One says
that \textit{$V'$ and $V''$ are twist equivalent}
\cite[Definition\,3.8]{inp-AndrSchn02}. Two rank $n$ \YD modules
$V'$ and $V''$ of diagonal type are called \textit{Weyl
equivalent}, if there exists a \YD module $V$ of diagonal type
which is twist equivalent to $V'$ as a braided vector space, and
an ordered basis $E''$ of $\ndZ ^n$ such that $E''\in E(V'')$ and
$V$ has degree $E''$ with respect to $V''$. Both twist equivalence
and Weyl equivalence are equivalence relations.
\end{defin}

\begin{bems}
1. The above definition doesn't depend on the ordering of $E''$
since the property of $V$ having degree $E''$ with respect to
$V''$ doesn't depend on it.

2. Weyl equivalence does not use the \YD structure, it is an
equivalence between braided vector spaces. In particular, \YD
modules, which are isomorphic as braided vector spaces, are always
Weyl equivalent.

3. In the theory of simple super Lie algebras it is a basic fact
that (in contrast to semi-simple Lie algebras) for different
choices of sets of simple roots one can obtain different Dynkin
diagrams. Weyl equivalence defined above reflects the aim to
consider objects which are independent of the choice of ``sets of
simple roots''.
\end{bems}

\begin{satz}
Let $V_1,V_2$ be Weyl equivalent \YD modules of diagonal type such
that both $\proots (\cB (V_1))$ and $\proots (\cB (V_2))$ are
finite. Then $\cB (V_1)$ and $\cB (V_2)$ have the same dimension
and the same Gel'fand--Kirillov dimension.
\end{satz}

\begin{bew}
If $V_1$ and $V$ are twist equivalent \YD modules then by
\cite[Proposition\,3.9]{inp-AndrSchn02} the corresponding Nichols
algebras are isomorphic as $\ndZ ^n$-graded vector spaces. In
particular, $\proots (\cB (V_1))=\proots (\cB (V))$. Further, if
$V$ has degree $E$ with respect to $V_2$ for some ordered basis
$E$ of $\ndZ ^n$ then by the first remark after
\cite[Proposition\,1]{a-Heck04c} one has $\roots (\cB (V))=\roots
(\cB (V_2))$ up to a change of basis in $\ndZ ^n$, and the heights
of the restricted \PBW generators coincide. This gives the claim.
\end{bew}

\section{Nichols algebras of rank 2}
\label{sec-rank2}

In \cite{a-Heck04a} finite dimensional rank 2 Nichols algebras
were described in terms of generators and relations using full
binary trees. In some of these examples one can drop the condition
on the order of some of the structure constants, without changing
$\proots (\cB (V))$. The (infinite dimensional) examples obtained
this way were already known in the literature.

Note that twist equivalent \YD modules of diagonal type are Weyl
equivalent. Therefore in order to analyze Weyl equivalence of rank
2 \YD modules of diagonal type one can assume that $q_{12}=1$. Let
$R_m$ denote the set of primitive $m^\mathrm{th}$ roots of unity.

\begin{satz}\label{s-rank2}
Let $V$ and $V'$ be \YD modules of diagonal type such that with
respect to certain bases their structure constants
$(q_{ij})_{i,j\in \{1,2\}}$ and $(q'_{ij})_{i,j\in \{1,2\}}$ are
listed in Figure 1. Then $V$ and $V'$ are Weyl equivalent if and
only if these constants appear in the same row with the same fixed
parameters.
\begin{figure}
 {\scriptsize
\begin{align*}
\begin{array}{l|l|l|c}
 \text{conditions for $q_{ij}$} & \text{free parameters} &
 \text{fixed parameters} & \text{tree}\\
\hline
 q_{21}=1 & & q_{11},q_{22}\in k\setminus \{0\} & T1\\
\hline
 q_{21}=q_{11}^{-1},\ q_{22}=q_{11} & & q_{11}\in k\setminus \{0,1\} & T2\\
\hline
 q_{11}=q,\ q_{21}=q^{-1},\ q_{22}=-1,\text{ or} & &
 q\in k\setminus \{-1,0,1\} & T2\\
 q_{11}=-1,\ q_{21}=q,\ q_{22}=-1 & & & T2\\
\hline
 q_{11}=q,\ q_{21}=q^{-2},\ q_{22}=q^2 & & q\in k\setminus \{-1,0,1\} & T3\\
\hline
 q_{11}=q,\ q_{21}=q^{-2},\ q_{22}=-1 & q\in \{q_0,-q_0^{-1}\} &
 q_0\in k\setminus \{-1,0,1\} & T3\\
\hline
 q_{11}=\zeta ,\ q_{21}=q^{-1},\ q_{22}=q & q\in \{q_0,\zeta q_0^{-1}\} &
 \zeta \in R_3 & T3\\
 & & q_0\in k\setminus \{0,1,\zeta ,\zeta ^2\} &\\
\hline
 q_{11}=\zeta ,\ q_{21}=-\zeta ,\ q_{22}=-1 & \zeta \in R_3 & & T3\\
\hline
 q_{11}=\zeta ^4,\ q_{21}=\zeta ^{-3},\ q_{22}=-\zeta ^2,\text{ or}&
 \zeta \in \{\zeta _0,-\zeta _0^{-1}\}&
 \zeta _0\in R_{12} & T4\\
 q_{11}=\zeta ^4,\ q_{21}=\zeta ^{-1},\ q_{22}=-1,\text{ or}& & & T5\\
 q_{11}=\zeta ^{-3},\ q_{21}=\zeta ,\ q_{22}=-1 & & & T7\\
\hline
 q_{11}=-\zeta ^2,\ q_{21}=\zeta ,\ q_{22}=-\zeta ^2,\text{ or}&
& \zeta \in R_{12} & T4\\
 q_{11}=-\zeta ^2,\ q_{21}=\zeta ^3,\ q_{22}=-1,\text{ or}& & & T5\\
 q_{11}=-\zeta ^{-1},\ q_{21}=\zeta ^{-3},\ q_{22}=-1 & & & T7\\
\hline
 q_{11}=\zeta ,\ q_{21}=\zeta ^{-2},\ q_{22}=-\zeta ^3,\text{ or}&
& \zeta \in R_{18} & T6\\
 q_{11}=-\zeta ^2,\ q_{21}=-\zeta ,\ q_{22}=-1,\text{ or} & & & T14\\
 q_{11}=-\zeta ^3,\ q_{21}=-\zeta ^{-1},\ q_{22}=-1 & & & T9\\
\hline
 q_{11}=q,\ q_{21}=q^{-3},\ q_{22}=q^3 & & q\in k\setminus \{-1,0,1\} & T8\\
 & & q\notin R_3 & \\
\hline
 q_{11}=\zeta ^2,\ q_{21}=\zeta ,\ q_{22}=\zeta ^{-1},\text{ or}&
& \zeta \in R_8 & T8\\
 q_{11}=\zeta ^2,\ q_{21}=-\zeta ^{-1},\ q_{22}=-1,\text{ or} & & & T8\\
 q_{11}=\zeta ,\ q_{21}=-\zeta ,\ q_{22}=-1 & & & T8\\
\hline
 q_{11}=\zeta ^6,\ q_{21}=-\zeta ^{-1},\ q_{22}=\zeta ^8,\text{ or} & &
\zeta \in R_{24} & T10\\
 q_{11}=\zeta ^6,\ q_{21}=\zeta ,\ q_{22}=\zeta ^{-1},\text{ or}& & & T13\\
 q_{11}=\zeta ^8,\ q_{21}=\zeta ^5,\ q_{22}=-1,\text{ or} & & & T17\\
 q_{11}=\zeta ,\ q_{21}=\zeta ^{-5},\ q_{22}=-1 & & & T21\\
\hline
 q_{11}=\zeta ,\ q_{21}=\zeta ^{-3},\ q_{22}=-1,\text{ or} &
\zeta \in \{\zeta _0,\zeta _0^{11}\}& \zeta _0\in R_5\cup R_{20} & T11\\
 q_{11}=-\zeta ^{-2},\ q_{21}=\zeta ^3,\ q_{22}=-1 & & & T16\\
\hline
 q_{11}=\zeta ,\ q_{21}=\zeta ^{-3},\ q_{22}=-\zeta ^5,\text{ or} & &
\zeta \in R_{30} & T12\\
 q_{11}=-\zeta ^3,\ q_{21}=-\zeta ^4,\ q_{22}=-\zeta ^{-4},\text{ or}& & & T15\\
 q_{11}=-\zeta ^5,\ q_{21}=-\zeta ^{-2},\ q_{22}=-1,\text{ or} & & & T18\\
 q_{11}=-\zeta ^3,\ q_{21}=-\zeta ^2,\ q_{22}=-1 & & & T20\\
\hline
 q_{11}=\zeta ,\ q_{21}=\zeta ^{-3},\ q_{22}=-1,\text{ or} & &
\zeta \in R_{14} & T19\\
 q_{11}=-\zeta ^{-2},\ q_{21}=\zeta ^3,\ q_{22}=-1 & & & T22\\
\hline
\end{array}
\end{align*}
}
\caption{Weyl equivalence for rank 2 \YD modules}
\end{figure}
\end{satz}

\begin{bew}
One has to show that if $V$ is a \YD module such that its
structure constants are listed in Figure 1, then for $i\in
\{1,2\}$ the \YD module $V'$ having degree $s_i(E_0)$ with respect
to $V$ is twist equivalent to a \YD module in the same row of
Figure 1. In this case let's say that $V'$ is connected with $V$
via $s_i$. Note that if $V'$ is connected with $V$ via $s_i$ then
also $V$ is connected with $V'$ via $s_i$. Further, it has to be
shown that for any two \YD modules $V,V'$ in the same row there is
a set $V=V_1,V_2,\ldots ,V_m=V'$ of \YD modules such that $m\ge 2$
and $V_i$ is connected with $V_{i+1}$ via some $s_{j_i}$, where
$j_i\in \{1,2\}$, whenever $1\le i<m$. The relation between the
structure constants of $V'$ and $V$ is given by Equations
(\ref{eq-qtrafo}). The calculations are straightforward. Note also
that sometimes one has to perform an additional base change
$\{x_1,x_2\}\to \{x_2,x_1\}$ of $V'$.

As an example consider the second last row of Figure 1. Let $V_i$,
where $1\le i\le 4$, denote the \YD module corresponding to the
$i^\mathrm{th}$ collection of structure constants. Then $V_1$ is
connected with $V_1$ via $s_1$, $V_1$ is connected with $V_3$ via
$s_2$ (perform base change $\{x_1,x_2\}\to \{x_2,x_1\}$ for
$V_3$), $V_3$ is connected with $V_4$ via $s_2$, $V_4$ is
connected with $V_2$ via $s_1$, and $V_2$ is connected with $V_2$
via $s_2$.
\end{bew}
\enlargethispage{\baselineskip }

\bibliographystyle{mybib}
\bibliography{quantum}

\end{document}